\documentclass[12pt,a4paper]{article}

\usepackage{amssymb,amsmath,amsthm,amsfonts,graphicx,color,epsfig}

\newtheorem{thm}{Theorem}[section]
\newtheorem{lemma}[thm]{Lemma}

\newtheorem{prop}[thm]{Proposition}

\theoremstyle{remark}

\newtheorem{remark}[thm]{Remark}

\newtheorem{definition}[thm]{Definition}

\newtheorem{example}[thm]{Example}

\def\R{{\mathbb R}}
\def\Z{{\mathbb Z}}

\def\Q{{\mathbb Q}}

\title{Makanin-Razborov diagrams for limit groups}
\author{Emina Alibegovi\'{c}}

\begin{document}
\renewcommand{\thefootnote}{\null}
\maketitle
\setcounter{footnote}{0}
\renewcommand{\thefootnote}{\arabic{footnote}}

\begin{abstract}

We give a description of $Hom(G,L)$, where $L$ is a limit group
(fully residually free group). We construct a finite diagram of groups,
Makanin-Razborov diagram, that gives a convinient representation of
all such homomorphisms. 

\end{abstract}

\section{Introduction}

The subject of this paper has the roots in the following problem:
given a finitely presented group $G$, describe the set of all
homomorphisms $G\to F$ to a fixed free group $F$.

When $G$ is the fundamental group of a closed surface,
say of genus $g$, Stallings answered our question as follows. Denote
by $q:G \to F_g$ an epimorphism to the free group of rank $g$ (defined
by inclusion of the boundary to the handlebody). Then every $f:G\to F$
factors as $f=\phi \circ q \circ \alpha$, for some automorphism
$\alpha: G \to G$, and some $\phi:F_g \to F$. Thus $Hom(G,F)$ is
`parametrized' by the product of the Teichm\"{u}ller modular group of
$G$ and the `affine space' $F^g$. This theorem of Stallings was
generalized to arbitrary finitely generated groups $G$ by Sela in
\cite{zlil1} and Kharlampovich and Myasnikov in \cite{km1}.

\begin{figure}[!ht]
\centerline{\input{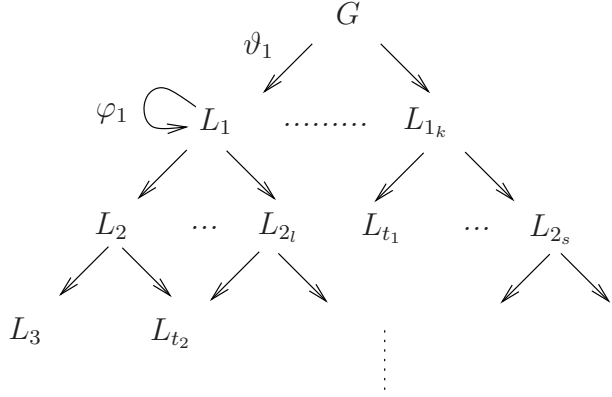}}
\caption{M-R diagram}\label{mr}
\end{figure}

All homomorphisms $G \to F$ are encoded into a finite diagram of
groups, called the {\it Makanin-Razborov diagram}. Each group in
this diagram has a finite number of directed edges issuing from
it. This number will be zero if the group in question is a free
group. Each edge represents a quotient map, and all quotients are
proper, see Figure \ref{mr}.

\noindent
Every homomorphism $h: G \to F$ can be written as:
\[h= h_0 \circ  \vartheta_{k} \circ \varphi_{k-1}
\circ \ldots \circ \vartheta_1 \circ \varphi_0, \] 
\noindent
where $\varphi_i \in Mod(L_i)$, $\vartheta_i$ are the quotient maps
and $h_0:L_{t_i}\to F$, where $L_{t_i}$ is a free group.  We say that
$h$ factors through a branch of the M-R diagram. In addition, all
groups in the M-R diagram naturally belong to the class of limit
groups. Limit groups, also known as fully residually free groups, have
been studied in \cite{zlil1}, \cite{km1}, \cite{km2},
\cite{lambdatrees}. The structure of limit groups can be easily
described. These groups can be built inductively: level 0 limit groups
are finitely generated free groups, finitely generated free abelian
groups and surface groups. A level $n$ limit group is obtained by
taking a finite number of free products or amalgamated free products
or HNN extensions of level $n-1$ limit groups along their maximal
cyclic subgroups. We will talk more about this decomposition in
Section \ref{back}. 

We are interested in studying $Hom(G,L)$, when $G$ is a f.g. group,
and $L$ is an arbitrary limit group, and constructing M-R diagrams for
this case. The problem that occurs here lies in the fact that a
homomorphism $h:G \to L$ might not factor through a complete branch of
the M-R diagram. That is, the final homomorphism $h_0$ in the above
representation of $h$ might be an embedding of some limit group into
$L$. We need to ensure that there are only finitely many such
embeddings, up to some equivalence relation. The following theorem is
the objective of this work:

\begin{thm}\label{mainthm}[Main Theorem]

Let $G$ be a f. g. freely indecomposable group and $L$ a
freely indecomposable limit group. There exist finitely many proper
quotients $G_1, \ldots, G_r$ of $G$ so that for every homomorphism
$f:G \to L$ with $d_{[f]}\gg 0$ an element of the equivalence class $\sim$ of
$f$ factors through some $G_i$. Furthermore, there are only finitely
many homomorphisms with uniformly bounded $d_f$ and nonabelian
image, up to conjugacy.
\end{thm}

We will define all the terms used in the statement of this theorem in
the next section, and the proof will be given in Section \ref{main}.

\section{Background}\label{back}

We list some of the properties of limit groups that are often
used, without proofs, and we refer reader to \cite{zlil1}, and \cite{mnotes}.

\begin{lemma}
Let $L$ be a limit group. 

\begin{itemize}
\item[(L0)] $L$ is torsion free.
\vspace{-.3cm} 
\item[(L1)] $L$ is finitely presented, in fact coherent. 
\vspace{-.3cm} 
\item[(L2)] Every f.g. subgroup of $L$ is a limit group. 
\vspace{-.3cm} 
\item[(L3)] Every abelian subgroup of $L$ is finitely generated and
  free, and there is a uniform bound on its rank. 
\vspace{-.3cm} 
\item[(L4)] Every abelian subgroup is contained in a unique maximal
  abelian subgroup.
\vspace{-.3cm} 
\item[(L5)] Every maximal abelian subgroup of $L$ is malnormal. 
\end{itemize}

\end{lemma}

In \cite{relhyp} we defined a class of groups that contains limit
groups, and for this class we found $\delta$-hyperbolic spaces on
which they act freely, by isometries. This class of groups was
defined as follows:

\begin{definition}\label{cgroup}
A torsion-free, f.g. group $G$ is a depth 0 $\mathcal
C$-group if it is either an f.g. free group, or an f.g. free abelian
group or the fundamental group of a closed hyperbolic surface. A torsion-free
f.g. group $G$ is a $\mathcal C$-group of depth $\leq n$ if it has a graph of
groups decomposition with three types of vertices: abelian, surface or
depth $\leq (n-1)$, cyclic edge stabilizers and the following holds:

\begin{itemize}
\item Every edge is adjacent to at most one abelian vertex $v$.
  Further, $G_v$, the stabilizer of $v$, is a maximal abelian subgroup
  of $G$.
\vspace{-.3cm}
\item Each surface vertex group is the fundamental group of a surface with
  boundary, and to each boundary component corresponds an edge of this
  decomposition. Each edge group is conjugate to a boundary
  component. 
\vspace{-.3cm}
\item The stabilizer of a depth $\leq (n-1)$ vertex $v$, $G_v$, is
  $\mathcal C$-group of depth $\leq (n-1)$. The images in $G_v$ of
  incident edge groups are distinct maximal abelian subgroups of $G_v$
  (i.e., cyclic subgroups generated by distinct, primitive, hyperbolic
  elements of $G_v$).
\end{itemize}

We say that the depth of a $C$-group $G$ is the smallest $n$ for which
$G$ is of depth $\leq n$. 
\end{definition}

That limit groups belong to the class $\mathcal{C}$ follows from
Theorem 3.2. and Theorem 4.1. in \cite{zlil1}. In fact,
the decomposition of a limit group from this definition coincides with
cyclic JSJ decomposition defined in \cite{zlil1}. 

Let $L$ be a depth $n$ limit group. Let $\Delta_L$ be a graph of
groups decomposition of $L$ given in the above definition, call $T_L$
the underlying graph, and let $T$ be the tree so that $T/L=T_L$. In
\cite{relhyp} we showed that we can find a graph of spaces $X/L$
corresponding to $\Delta_L$ so that its universal cover $X$ is
$\delta$-hyperbolic. In order to establish necessary notation we give
some properties of the space $X$. $X/L$ is quasiisometric to the wedge
of $k$ rays $[0,\infty)$ joined at 0. We lift $k$ rays that correspond
to $\partial X/L$ to rays $r_i:[0,\infty) \to X, \ i=1,\ldots,k$, and
we let $h_i$ be the horofunction corresponding to $r_i$.  The
stabilizer, $L_i < L$, of the limit point $r_i(\infty)$ of $r_i$
preserves $h_i$. Denote by $B_i(\rho)$ the horoballs
$h_i^{-1}(-\infty,\rho)\subset X$. For sufficiently small $\rho$ the
intersection $\gamma B_i(\rho) \cap B_j(\rho)$ is empty unless $i=j$
and $\gamma \in L_i$. Let
\[LB(\rho)=\bigcup_{i, \gamma}\gamma B_i(\rho), \] 
$i=1, \ldots, k, \ \gamma \in L$. Let $X(\rho)=X\backslash
LB(\rho)$. $X(\rho)/L$ is compact for all $\rho \in (-\infty,
\infty)$.

These properties in fact constitute the definition of relatively
hyperbolic groups given by Gromov in \cite{gromov}. Hence limit groups
are hyperbolic relative the collection of the representatives of
conjugacy classes of their maximal noncyclic abelian subgroups (see
\cite{relhyp}). We will call the subgroups $L_i$ parabolic subgroups.

The {\it modular group} $Mod(L)$ associated to the
decomposition $\Delta_L$ of $L$ is the subgroup of $Aut(L)$ generated by
\begin{itemize}
\item inner automorphisms of $L$,
\vspace{-.3cm}
\item Dehn twists in the centralizers of edge groups,
\vspace{-.3cm}
\item automorphisms induced by automorphisms of abelian vertex groups
that are identity on peripheral subgroups and all other vertex groups,
and
\vspace{-.3cm}
\item automorphisms induced by homeomorphisms of surfaces underlying
surface vertex groups that fix all boundary components.
\end{itemize}

Let $G$ be an f.g. group with a finite generating set $S$. We consider
a sequence of homomorphisms $f_i:G \to L$. Each of the given
homomorphisms induces an action of $G$ on the space $X$, and for each
we define

$$d_i=\inf \{ \sup \{d(x,f_i(g)x): g \in S \}:\ x \in X\}.$$ 
\noindent
If this infimum is attained at a point $x_i$ then $x_i$ is called a
centrally located point for the action.

In order to apply \textit{Compactness Theorem}, \cite{rtrees}, i.e.,
to see if we can extract a subsequence of actions that converges, we
need to see whether each of these actions has a centrally located
point and how the sequence ${d_i}$ behaves.

The following lemma has been proved by Bestvina in \cite{degen} for 
hyperbolic spaces $\mathbb{H}^n$ and by Paulin in \cite{paulin} for 
Gromov hyperbolic spaces. 
 
\begin{lemma}\label{clp}
For every $f:G \to L$ whose image is not an abelian group and for
which an induced action on $X$ is nonelementary (no point at infinity
is fixed by the whole group), there exists a centrally located point.

\end{lemma}

\begin{proof}
Suppose $d_f=0$. Since the group $L$ contains no elliptic elements,
this implies that for every $g \in S$, $f(g)$ is a parabolic
element. Furthermore, $f(G)$ is an abelian subgroup of $L$. Note that
a partial converse holds: if $f(G)$ is an abelian group and all
generators are mapped into parabolic elements then $d_f=0$ and
centrally located point does not exist.  So we need to show that if
$d_f>0$ then a centrally located point exists.

We would like to show that a map $F:X\to \R_+$ defined by 
$$F(x)=\sup_{g \in S} d(x,f(g)x)$$ is a proper map away from the horoballs
and consequently attains its infimum. Suppose it is not, that is there
is a sequence $\{x_n\}$ in $X$ not contained in a compact set, such
that $\{F(x_n)\}$ is bounded. Recall from the definition of $X$ that
there is $\rho$ such that 
$B_i(\rho) \cap \gamma B_j(\rho)= \emptyset$,  
unless $i=j$ and $\gamma \in L_i$ and $X(\rho)/L$ is compact, where 
$X(\rho)=X\backslash \bigcup_{i,\gamma}\gamma B_i(\rho).$
We first note that our sequence has to
stay within bounded distance from $X(\rho)$, for
otherwise $F(x_n)\to \infty$ since $d_f>0$. Namely, suppose
$d(x_n,X(\rho))\longrightarrow \infty$ as $n \to \infty$. We also may
assume that $x_n \in B_i(\rho)$. Since $d_f >0$ not every $f(g), g \in
S$, is contained in $L_i$. Let us assume that $h \in S$ is such that $f(h)
\notin L_i$. Then $f(h)B_i(\rho)\cap B_i(\rho)= \emptyset$  implies
that $d(x_n, f(h)x_n) \geq 2d(x_n,X(\rho))$, hence $F(x_n) \rightarrow
\infty$ as $n\to \infty$. Contradiction. After passing to a
subsequence, $\{x_n\}$ will converge to a point $x \in \partial
X$. Since $d(x_n,f(g)x_n)$ is bounded for every $g \in S$ , we conclude
that the sequence $\{f(g)x_n\}$ also converges to $x$. Hence, $x$ is
fixed by $f(G)$, which is a contradiction.
\end{proof}

\begin{prop}\label{finite}
There are only finitely many homomorphisms $f:G \to L$ with
(uniformly) bounded $d_f$ and nonabelian image, up to conjugation. 
\end{prop}

\begin{proof}

Suppose this is not true, and there are infinitely many homomorphisms
$f_i:G \to L$ with $d_i$ bounded by $D >0$. We assume that
no $d_i = 0$, since otherwise $f_i$ has an abelian image, by Lemma
\ref{clp}. By the same lemma we know that a centrally located point
$x_i$ exists for the action given by $f_i$. 

To simplify the notation we will assume that $X/L$ has only one cusp.
Let $r:[0,\infty) \to X$ be a ray in $X$ corresponding to this cusp,
and let $h$ be the horofunction corresponding to $r$. Also, let $A$ be
an abelian subgroup of $L$ that stabilizes the horoball
$h^{-1}(-\infty,0)$. Choose $\rho$ small enough so that $B(\rho)\cap
\gamma B(\rho)= \emptyset$, for all $\gamma \in L\backslash A$. We
first note that $x_i$ can not be too deep inside the horoball. By
'deep` we mean that the distance from $x_i$ to the boundary of the
horoball has to be smaller than $D$. If it is not, then the ball of
radius $D$ around $x_i$ is not only completely contained within the
horoball, but also contains $f_i(g)x_i$, for all $g \in S$, since $d_i
\leq D$. This implies that $f_i(G)$ is abelian. Hence, the $D$-ball
around $x_i$, call it $B_i$, has to intersect $X(\rho)$. We consider
$X(\rho -D)$. The action of $L$ on $X(\rho-D)$ is cocompact, hence we
can find a compact set $K$ whose translates cover $X(\rho-D)$, see Figure
\ref{fundomain}. We can find $l_i \in L$  for each $x_i$ so that $l_ix_i
\in K$.  Since $K$ is compact there exist $r>0$ and $x \in K$ so that
$B_r(x)$ contains the translates $l_iB_i$, for all $i$.
\begin{figure}[!ht]
\centerline{\input{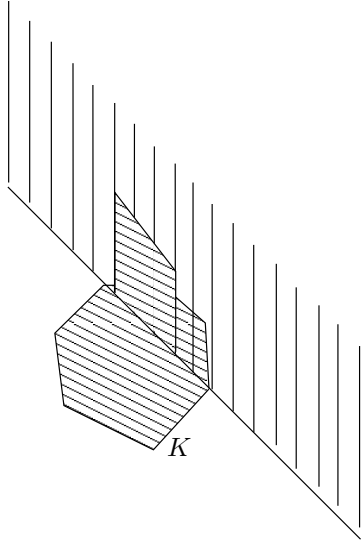}}
\caption{Fundamental domain for the action of $L$}\label{fundomain}
\end{figure}

\begin{eqnarray*}
     d(x,l_if_i(g)l_i^{-1 }x) & = & d(l_i^{-1}x,f_i(g)l_i^{-1}x) \leq \\ 
& \leq & d(l_i^{-1}x,x_i)+d(x_i,f_i(g)x_i)
+d(f_i(g)x_i,f_i(g)l_i^{-1}x) \leq \\ 
& \leq & 2r+D
\end{eqnarray*} 

\noindent
and so $l_if_i(g)l_i^{-1}$ moves a point $x$ within a ball of radius
$2r+D$, for all $i, \ \text{and for all} \ g \in S$. Since there are
only finitely many translates of $x$ within that ball we conclude that
there can be only finitely many nonconjugate $l_if_il_i^{-1}$.
\end{proof}

We therefore conclude that if homomorphisms $f_i$ belong to distinct
conjugacy classes and have nonabelian images, the sequence of actions
they induce contains a subsequence which converges to an action of $G$
on an $\R$-tree $T_{\infty}$ without a global fixed point. Let us call this
limiting action $\rho$. 

\begin{figure}[!ht]
\centerline{\input{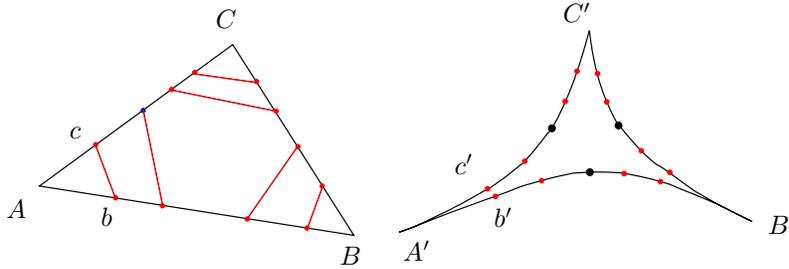}}
\caption{Defining the lamination $\widetilde{\Lambda_f}$ on
$\widetilde{K}$}\label{lam}
\end{figure}

For every homomorphism $f:G \to L$ we get a measured lamination
$\Lambda_f$ on the complex $K$ whose fundamental group is $G$. We
define a resolution $\Phi: \widetilde{K} \to X$ by defining it $f$ -
equivariantly on the vertices of the triangulation of
$\widetilde{K}$. We extend it to edges so that each edge of the
triangulation is mapped into the geodesic between the images of its
endpoints. Finally extend the map equivariantly to 2-cells.  The image
of each triangle has a unique measured lamination on it. We will pull back this
lamination to get $\widetilde{\Lambda_f}$, see Figure \ref{lam}. Say $ABC$
is a triangle in $\widetilde{K}$ whose image under $\Phi$ is a
triangle $A'B'C'$ in $X$. Suppose $b'$ and $c'$ are the points on
$A'B'$ and $A'C'$, respectively, the same distance from $A'$.

\noindent
Let $b$ and $c$ be the points on $AB$ and $AC$,
respectively, for which $\Phi(b)=b'$ and $\Phi(c)=c'$, and let $[bc]$
be the segment for which $\Phi([bc])=[b'c']$. The segment $[bc]$ is
contained in a leaf of our lamination $\widetilde{\Lambda_f}$. Notice
that the length of the image of $[bc]$ is smaller than or equal to
$6\delta$. $\Lambda_f$ is then a projection of $\widetilde{\Lambda_f}$
to $K$. 

\begin{definition}
For a measured lamination $(\Lambda,\mu)$ on $K$ we define
the {\it length} of $\Lambda$ to be
\[L(\Lambda)=\sum_{e \in K^{(1)}}\mu(e).\]
The length of a resolution $\Phi:\tilde{K}\to X$, $L(\Phi)$, is the
length of the induced lamination.
\end{definition}

For each homomorphism $f:G \to L$ we choose an $(f)$-resolution
$\Phi_f$ so that

$$L(\Phi_f) = \inf \{L(\Phi): \Phi \ \text{is an}\ (f)-
\text{resolution}\}. $$
We need to verify that such a resolution $\Phi_f$ exists. 

\begin{lemma}
Let $\{\Phi_i\}$ be a sequence of resolutions with
bounded lengths for a homomorphism $f:G \to L$ with $d_f >0$. There
exists a resolution $\Phi: \widetilde{K} \to X$ such that
$L(\Phi)=\lim L(\Phi_i)$. 
\end{lemma}

\begin{proof}
Let $K'$ be the fundamental domain for the
action of $G$ on $\widetilde{K}$. We first note that $\lim L(\Phi_i)$
has to be positive. If it were not, we would have that $d(X(\rho),
\Phi_i(K')) \to \infty$ as $i \to \infty$. On the other hand, since
the image of $f$ is not abelian, there has to exist $g \in G$ so that
$\Phi_i(gK')$ belongs to a horoball disjoint from the one containing
$\Phi_i(K')$. Hence, there should also exist a translate of
$\Phi_i(K')$ that intersects $X(\rho)$ nontrivially, since
$\Phi_i(\widetilde{K})$ is connected, but that is impossible under
these hypotheses. Therefore, $L=\lim L(\Phi_i) > 0$. By definition
\[L(\Phi_i)=\sum_{j=1}^{E} \mu_i(e_j), \]
where $E$ is the number of edges in $K^{(1)}$, and $(\Lambda_i,
\mu_i)$ is a measured lamination induced by $\Phi_i$. After passing to a
subsequence, we may assume that for every $j=1, \ldots, E$, the
sequence $\{\mu_i(e_j)\}$ converges, and we denote its limit by
$\ell_j$. Hence, 
\[L=\sum_{j=1}^E \ell_j. \]
Since $\{L_i\}$ is a bounded sequence and its limit is positive, we
know that all of the vertices of $K'$ are mapped into a bounded
neighborhood of $X(\rho)$ by each $\Phi_i$. Pick one of these
vertices, say $v$. We can assume that $v$ is contained in $X(\rho)$
(if it is not, we can always enlarge $X(\rho)$ so as to cover
$\Phi_i(v)$). If we denote by $X'$ the fundamental domain of the
action of $L$ on $X(\rho)$, then there exists $l_i\in L$ such that
$l_i\Phi_i(v) \in X'$. Furthermore, $l_i\Phi_i(K')$ is contained in a
bounded neighborhood of $X'$. Hence, for every vertex $v$ of $K'$ we
can find a convergent subsequence of $l_i\Phi_i(v)$. We define
\[\Phi(v)=\lim_{i \to \infty} l_i\Phi_i(v), \ v \ \text{vertex of}\ \ 
K'.\]
We extend $\Phi$ $f$-equivariantly, and get a resolution. Clearly,
$\Phi$ has the desired length, since $l_i$'s are isometries and
preserve the lengths of $\Phi_i$'s. 

\end{proof}

We are interested in the closure, $LIM(K)$, of all $\Lambda_f$'s
inside the space of projectivized measured laminations,
$\mathcal{PML}(K)$. The space $LIM(K)$ is compact.

\begin{remark}
Let us back up to the sequence $\{f_i\}$ which converges to an action
$\rho$ of $G$ on a tree $T_{\infty}$. For each of $f_i$'s we take a short
resolution $\Phi_{f}$ and a corresponding lamination
$\Lambda_{f_i}$. This sequence will converge in $\mathcal{PML}(K)$
to some $\Lambda$. There exists a resolution $\Phi:\widetilde{K} \to
T_{\infty}$ so that the leaves of $\Lambda$ are mapped to points. 
\end{remark}
\begin{definition}
Let $f:G \to L$ be a homomorphism between two limit groups. We define
a collection of moves that are allowed to be performed on $f$:

\begin{enumerate}

\item[(M1)] precompose $f$ by an element of $Mod(G)$, 
\vspace{-.2cm}
\item[(M2)] conjugate $f$ by an element of $L$, 
\vspace{-.2cm}
\item[(M3)] if there is an abelian vertex group $A$ in $\Delta_G$ such that
  $f(A)$ is contained in a parabolic subgroup of $L$, then redefine $f$ on
  $A$ so that the new homomorphism coincides with $f$ on the adjacent edge
  groups,
 \vspace{-.2cm}
\item[(M4)] bending: suppose there is an edge group $E$ in $\Delta_G$
  whose image is nontrivial and contained in a parabolic subgroup of
  $L$. If the edge $X_E$ corresponding to $E$ separates $\Delta_G$, we
  conjugate the image of one of the connected components of $\Delta_G
  \backslash X_E$ by an element of $L$ that commutes with $f(E)$. If
  edge $X_E$ is nonseparating, i.e., corresponds to an HNN extension,
  we multiply the image of the Bass-Serre generator by an element of
  $L$ that commutes with $f(E)$.
  
\end{enumerate}

\noindent
We say that homomorphisms $f, g:G \to L$ are equivalent, $f \sim g$, if
there is a sequence  $f=f_0, f_1, \ldots, f_n=g:G \to L$
such that $f_{i+1}$ is obtained from $f_i$ by performing one of the moves
(M1)-(M4). 

\end{definition}

\begin{remark}
If we take $L$ to be a free group then the only
equivalent homomorphisms are the ones that either differ by an element of
$Mod(G)$ or are conjugates of each other. Why did we need to add moves
(M3) and (M4)? Let us look at the following example. 
\end{remark}

\begin{example}\label{badex}
Let $G=<a,b,s,t|w(a,b)=s, [s,t]=1> =F_2*_{\Z}\Z^2$, where the word
$w(a,b)$ is chosen so that $G$ is a limit group. Consider the sequence
of homomorphisms $f_i:G \to G$ defined by $f_i(a)=a, f_i(b)=b,
f_i(s)=s$, and $f_i(t)=t^i$. This sequence has the property that
$d_i\to \infty$, but all of its members are embeddings. This is where
we will need move (M4).  
\end{example}

\begin{definition}
A homomorphism $f$ is {\it short} if $\Phi_f$ is shortest
among all $\Phi_g$ when $g \sim f$. $LIM'(K)$ will denote the
closure of the set of $\Lambda_f$'s, for short $f$'s, in
$\mathcal{PML}(K)$. 
\end{definition}
That the shortest homomorphism exists can be
proved in the exactly same way as Proposition \ref{finite}. 

\section{Proof of Main Theorem}\label{main}

Let us first note that it is sufficient to consider only groups which
are $\omega$-residually free groups. If $G$ is not $\omega$-residually
free, then there are nontrivial elements $g_1, \ldots, g_k$ so that
every $f:G \to F$ kills at least one. We claim that every homomorphism
$G \to L$ kills at least one of the $g_i$'s. Suppose not, i.e., there
is a homomorphism $\phi:G \to L$ such that $\phi(g_i)\neq 1$, for all
$i=1, \ldots,k$. Since $L$ is $\omega$-residually free group, we can
find a homomorphism $\rho: L \to F$ such that $\rho(\phi(g_i))\neq
1$. We have found a homomorphism $\rho \phi:G \to F$ which is
injective on the set $\{g_1, \ldots, g_k\}$, a contradiction. Hence
the quotients $G\to G/\ll g_i\gg$ satisfy the requirements of Main Theorem.

From now on we fix a complex $K$ that reflects the decomposition of the limit
group $G$ as in Definition \ref{cgroup}. We will prove our theorem by
considering the space $LIM'(K)$. Our first concern is to find suitable
proper quotients of $G$. The series of lemmas that follow will show
how to obtain these quotients depending on the type of the individual
laminations in $LIM'(K)$. We first concentrate on a lamination $\Lambda$
which is a limit of laminations $\Lambda_{f_i}\in LIM'(K)$, where $f_i:G\to L$
are homomorphisms with nonabelian images and the property that $d_i
\to \infty$. 

\begin{lemma}\label{simplicial}
If $\Lambda$ is entirely simplicial we form the following quotients: 
\begin{itemize}
\vspace{-.2cm}
 \item Abelianize the subgroup carried by a generic leaf in each
    Cantor set bundle, if these subgroups are nonabelian.
 \vspace{-.2cm}   
 \item If all of these subgroups are abelian, for each of them we make the
    following quotients: 
\vspace{-.2cm}       
\begin{itemize} 
\vspace{-.1cm} 
        \item mod out by its normal closure, 
\vspace{-.1cm}          
\item abelianize each vertex group in the splitting of $G$
           inherited from $\Lambda$.
  \end{itemize}
\end{itemize}
\vspace{-.2cm}
At least one of these quotients is proper. 
\end{lemma}

\begin{proof}
  We will assume that there is only one family of parallel leaves and
  let $N$ be its regular neighborhood.  $N$ is homeomorphic to $\ell
  \times [0,1]$, where $\ell$ is a leaf of the lamination $\Lambda$.
  Denote by $\ell_1$ and $\ell_2$ the boundary components of $N$, see
  Figure \ref{simpl}. Let $H$ be a subgroup of $G$ carried by
  $\ell_1$. $H$ cannot be trivial since that would imply that $G$ is
  freely decomposable.  Suppose $H=<h_1, \ldots, h_k>$.

\vspace{.2cm}
{\it (i)} If $H$ is nonabelian, so are its conjugates, hence
abelianization of a subgroup carried by a generic leaf yields a
proper quotient.    

\vspace{.2cm}
{\it (ii)} If $H$ is abelian and $f_i(H)=1$, then our homomorphism
factors through a proper quotient $G/\ll H \gg$. Hence we may assume
that there is at least one $j \in \{1, \ldots, k\}$ for which
$f_i(h_j)\neq 1$. Since $H$ is an abelian group, $f_i(H)$ is contained
either in a parabolic subgroup of $L$ of rank $n$ or in a cyclic
subgroup of $L$ generated by a hyperbolic element. 
 
\begin{figure}[!ht]
\centerline{\input{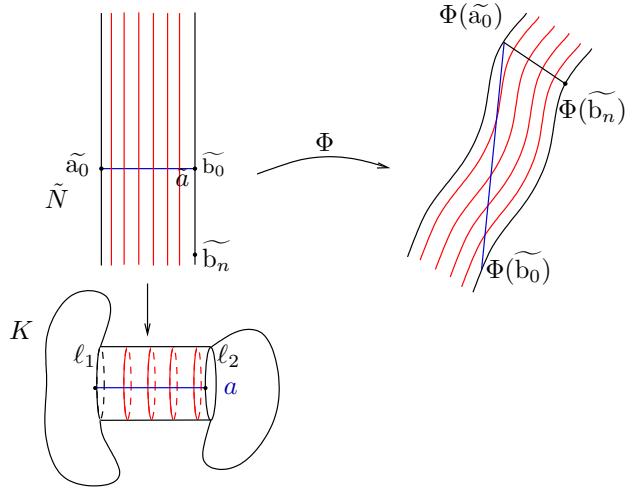}}
\caption{Shortening the simplicial component}\label{simpl}
\end{figure} 
\noindent

We would like to know what the general position of
    $\Phi_i(\widetilde{N})$ is within $X$. Due to the equivariance,
    the first thing we conclude is that $\Phi_i(\widetilde{N})$ is
    contained in a bounded neighborhood of a horoball $B$ that
    corresponds to the parabolic subgroup of $L$ which contains the
    image of $H$. Further, we can also conclude that, since
    $\mu_i(\ell_1)/d_i \to 0$, all the points of
    $\Phi_i(\widetilde{\ell}_1)$ are approximately at the same
    distance from $P$ (same for $\Phi_i(\widetilde{\ell_2})$).

\vspace{.2cm}
{\bf {\small $\Phi_i(\widetilde{\ell}_1)$ and
$\Phi_i(\widetilde{\ell}_2)$ are completely contained in $B$}}

Suppose that $\Phi_i(\widetilde{\ell}_2)$ is at a distance
proportional to $d_i$ from $P$. Assume that the fundamental group of
the vertex space, call it $K_2$, which is adjacent to $\ell_2$, is not
abelian. If $f_i(\pi_1(K_2))$ is abelian then $f_i$ factors through
the proper quotient of $G$ obtained by abelianizing $\pi_1(K_2)$ and
our claim holds. Otherwise, we can find a closed path $p \subset K_2$
whose measure is close to 0 and for which $[f_i(p),f_i(\ell_2)] \neq
1$ (we use the same notation for the elements of the fundamental group
and loops they represent).  Consider the image under $\Phi_i$ of a
lift $\tilde{p}$ of $p$.  One of the endpoints of $\Phi_i(\tilde{p})$
lies on $\Phi_i(\widetilde{\ell}_2)$, but $\Phi_i(\tilde{p})
\nsubseteq B$, for if it were we could conclude that $f_i(\ell_2)$ and
$f_i(p)$ commute. Therefore the length of $\Phi_i(\tilde{p})$ is
proportional to $2d_i$. On the other hand, $\mu_i(p)$ is virtually
zero, contradiction.  Hence, if $f_i(\pi_1(K_2))$ is not abelian, the
distance from $\Phi_i(\widetilde{\ell}_2)$ to $P$ is much smaller than
$d_i$. The same holds for $\Phi_i(\widetilde{\ell}_1)$.
\noindent
We now have two cases to consider:
\begin{itemize}
\item[(1)] $f_i(\pi_1(K_j)), \ j=1,2$, is not abelian, and
\vspace{-.2cm}
\item[(2)] at least one of $f_i(\pi_1(K_j)), \ j=1,2$ is abelian.
\end{itemize}
\noindent
Our strategy is as follows. We consider an edge $a$ of the
  triangulation of $K$ that intersects leaves of the lamination
  $\Lambda$. We will actually show more than we claimed to be true:
  all the quotients that we make are proper and $f_i$ factors through
  at least one of them. If not, we will be able to decrease the measure of
  the edge $a$, hence obtaining a lamination shorter than $\Lambda_{f_i}$.

\vspace{.2cm}
\noindent
(1) Let $a$ be an edge in $K^{(1)}$ that intersects $N$ nontrivially,
and let $a \cap \ell_1=\{\text{a}_0\}$ and $a \cap
\ell_2=\{\text{b}_0\}$. We know now that $\Phi_i(\widetilde{\ell}_1)$
and $\Phi_i(\widetilde{\ell}_2)$ are contained in the
$n_i$-neighborhood of $P$, and $n_i/d_i \to 0$ as $i \to \infty$. If
we assume that the euclidean height of $P$ is 1, which we may without
loss of generality, then for $x,y \in N_{n_i}(P)$ we will have $|\ln
x_{n+1}-\ln y_{n+1}|\leq n_i$.

Let $\tau_i$ denote the maximal displacement of points on
$\Phi_i(\widetilde{\ell}_1)$ under the action of $f_i(\hat{H})$, where
$\hat{H}$ denotes the centralizer of $H$. We know that $\tau_i \leq
\mu_i(\ell_1)$ and that $d(\Phi_i(\widetilde{\text{a}}_0),
\Phi_i(\widetilde{\text{b}}_0)) \sim d_i$. If there were a point in
the $f_i(h_j)$ - orbit of $\Phi_i(\widetilde{\text{b}}_0)$ whose
distance to $\Phi_i(\widetilde{\text{a}}_0)$ is smaller than that of
$\Phi_i(\widetilde{\text{b}}_0)$, then we could precompose $f_i$ by an
appropriate Dehn twist determined by $h_j$ and obtain a resolution
shorter than $\Phi_i$, a contradiction (see Figure \ref{simpl}). Hence,

\[d(\Phi_i(\widetilde{\text{a}}_0), \Phi_i(\widetilde{\text{b}}_0))
\leq d(\Phi_i(\widetilde{\text{a}}_0), f_i(h_j)^m
\Phi_i(\widetilde{\text{b}}_0)), \, \forall m\in \Z.\] 

The orthogonal projections of $\Phi_i(\widetilde{\ell}_1)$ and
$\Phi_i(\widetilde{\ell}_2)$ to $P$  will be
contained in small Hausdorff neighborhoods of parallel real lines. There exists
an element $g$ of the parabolic subgroup of $L$ corresponding to $B$
which will translate these lines towards each other.  For $i$ large
enough, the translation length of $g$, call it $\tau_P(g)$, while
acting on $P$ is going to be `much' smaller than $d_i$.
We can find $m, k \in \Z$ and a point $y$ on $\Phi_i(\widetilde{\ell}_1)$
so that
\[d(g^m\Phi_i(\widetilde{\text{b}}_0), y') \leq \tau_P(g)\ \ \text{and}\ \ d(f_i(h_j)^ky,\Phi_i(\widetilde{\text{a}}_0)) \leq \tau_i. \]
where $y'$ is the orthogonal projection of $y$ onto the horizontal
hyperplane in $B$ containing $\Phi_i(\widetilde{\text{b}}_0)$, see Figure
\ref{bend}. We then have

\begin{figure}[!ht]
\centerline{\input{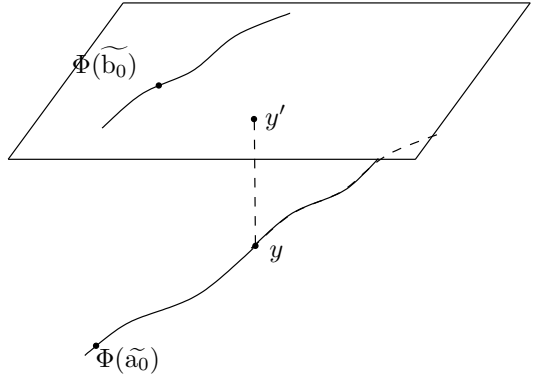}}
\caption{Shortening via bending}\label{bend}
\end{figure}

\begin{eqnarray*}
d(f_i(h_j)^kg^m\Phi_i(\widetilde{\text{b}}_0),\Phi_i(\widetilde{\text{a}}_0))
&\leq& d(g^m\Phi_i(\widetilde{\text{b}}_0),y')+d(f_i(h_j)^ky',\Phi_i(\widetilde{\text{a}}_0)) \leq\\ 
& \leq & \tau_P(g)+d(f_i(h_j)^ky',f_i(h_j)^ky)+ \\ 
& + & d(f_i(h_j)^ky,\Phi_i(\widetilde{\text{a}}_0)) \leq \\
& \leq & \tau_P(g)+\Big|\ln \frac{(\Phi_i(\widetilde{\text{a}}_0))_{n+1}}
 {(\Phi_i(\widetilde{\text{b}}_0))_ {n+1}}\Big| + \tau_i \ll d_i.
\end{eqnarray*}

\noindent
Therefore, performing the bending (a power of it, to be
more precise) determined by $g$ and precomposing the obtained
homomorphism with a power of the Dehn twist determined by $h_j$ will,
contrary to our assumption, give us a homomorphism in the same class
with $f_i$ which has a shorter resolution.

\vspace{.2cm} 
\noindent
(2) In the case that one of $f_i(\pi_1(K_j)), \, j=1,2$, is abelian we
may in fact assume that $\pi_1(K_j)$ is abelian itself. If not, $f_i$
will factor through the quotient of $G$ obtained by abelianizing
$\pi_1(K_j)$. Further, only one $K_j$ can have abelian fundamental
group, and without loss of generality we may assume it is $K_1$.  As
noted before $\Phi_i(\widetilde{\ell}_2)$ lies in the
$n_i$-neighborhood of $P$. If we further have that
$\Phi_i(\widetilde{\ell}_1)$ also lies in the $n_i$ neighborhood of
$P$, we apply (1). We therefore may assume that the distance from
$\Phi_i(\widetilde{\text{a}}_0)$ is proportional to $d_i$. 
$\Phi_i(\widetilde{\ell}_1)$ being mapped so deep into the horoball is
a consequence of $\mu_i(p)$, for every loop $p$ in $K_1$, being very short
compared to both $\mu_i(\ell_2)$ and $d_i$. We consider a new
$f_i$-resolution $\Phi_i':\widetilde{K}\to X$ which coincides with $\Phi_i$
on $\widetilde{\ell}_2$, but lowers $\widetilde{\ell}_1$. We explain the term 
`lowers' formally: there are constants $c_1, \ldots, c_n \in \R$ so
\def\Q{{\mathbb Q}}that for every point $x \in \widetilde{\ell}_1$
there is a point $x' \in \widetilde{\ell}_2$ such that
$(\Phi_i(x))_j=(\Phi_i(x'))_j + c_j, \ j=1, \ldots, n$. We define
$\Phi_i'(x)$ to be the point with coordinates 
\[((\Phi_i(x))_1, \ldots, (\Phi_i(x))_n, (\Phi_i(x'))_{n+1}).\]  
\noindent
Notice that we did not change the homomorphism $f_i$ in any way, we
only changed the resolution. By doing that we shortened the edge $a$
substantially. Unfortunately we increased the measures
of closed paths $p \subset K_1$ and $\ell_1$. However, if we
perform move (M3) on $f_i$ we obtain a homomorphisms with a shorter
resolution. Namely, the measure of the path $a$ will be the same as
under the new resolution for $f_i$, but the paths $p$ will now have
measures 0. This contradicts our assumption.

This concludes our discussion when
$\Phi_i(\widetilde{\ell}_1)$ and $\Phi_i(\widetilde{\ell}_2)$ are
contained in $B$. Further, these arguments remain valid when the
intersection of $\Phi_i(\widetilde{N})$ with the boundary horosphere
$P$ is a disjoint union of segments. 

\vspace{.2cm}
\noindent
{\bf {\small $\Phi_i(\widetilde{\ell}_j), \ j=1,2,$ are contained in $X(\rho)$ }} 
If $\Phi_i(\widetilde{\ell}_j), \ j=1,2,$ belong to different
connected components of $X(\rho)$, or if $\Phi_i(\widetilde{N}) \cap B
\neq \emptyset$ we can shorten the resolution by shortening, as above,
the distance between the connected components of
$\Phi_i(\widetilde{N}) \cap P$.  Suppose then that
$\Phi_i(\widetilde{N})$ is contained in a single connected component
of $X(\rho)$. $\Phi_i(\widetilde{N})/L$ is an annulus $F$ in a
compact part of $X/L$ with boundary components of length $\ll d_i$ and
of length approximately $d_i$. Without loss of generality we may
assume that $N$ had a drum triangulation, and $\Phi_i(\widetilde{N})$
has the same. Further we may assume that the triangulation on $N$ is
formed out of $k$ triangles. Each of these triangles has two sides
that are extremely long, meaning their lengths are greater than or
equal to $d_i-\mu_i(\ell_1)$, or $d_i-\mu_i(\ell_2)$, depending on
which boundary component of $N$ the short side of the triangle
lies. Let $s_0, \ldots, s_k$ denote the lifts of these long sides of
the triangles belonging to the same lift of $F$ in
$\Phi_i(\widetilde{N'})$, where $s_k$ is a translate of $s_0$ under an
element $g \in L$. Let $x$ be a point on $s_0$ at distance at least
$\mu_i(\ell_1)+k\delta$ and $\mu_i(\ell_2)+k\delta$, along $s_0$, from
$\Phi_i(\widetilde{\ell}_1)$ and $\Phi_i(\widetilde{\ell}_2)$,
respectively. There is a point $x_1$ on $s_1$ at a distance less than
or equal to $\delta$ from $x$ ($\delta$-hyperbolicity of
$X$). We continue moving through all the triangles in
$\Phi_i(\widetilde{N'})$, and conclude that there is a point $x_k \in
s_k$ at a distance less than or equal to $k\delta$ from $x$. If the
distance along $s_0$ between $x$ and $g^{-1}x_k$ is greater than
$k\delta$ then performing a Dehn twist will shorten the side $s_0$,
and all the others, and hence we obtain a shorter
resolution. 

\begin{figure}[!ht]
\centerline{\input{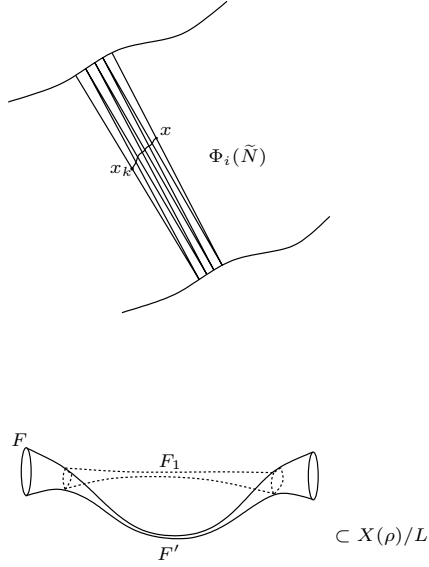}}
\caption{Shortening of the annulus $F$}\label{annuli}
\end{figure}
\noindent
Therefore, for the points $x, x_k$ as above
$d_{s_0}(x,g^{-1}x_k)<k\delta$. For every $j=0, \ldots,k$, let $s_j'$ denote
the projection to $F$ of a segment obtained from $s_j$ by removing the
ending subsegments of length $\mu_i(\ell_1)+k\delta$ and
$\mu_i(\ell_2)+k\delta$ that are adjacent to
$\Phi_i(\widetilde{\ell}_1)$ and $\Phi_i(\widetilde{\ell}_2)$,
respectively. For each $x \in s_0'$ we then have a loop based at $x$
of length $\leq 2k\delta$. Hence, boundary loops of $F$ are homotopic
to loops of length $\leq 2k\delta$, and there are only finitely many
such. Moreover, there is a constant $C$ so that for each two homotopic
loops of length $\leq 2k\delta$ there is a homotopy between them of
length $\leq C$. Let $\gamma_1$ and $\gamma_2$ be two of the short
loops we have found; the ones based at the initial and the end point
of $s_0'$, respectively. Also, let $F'$ denote the part of $F$ between
them. As we said earlier there is an annulus $F_1$ between $\gamma_1$
and $\gamma_2$ of length $\leq C$, given by the forementioned short
homotopy. We now replace $F'$ by $F_1$, see Figure \ref{annuli}. We
obtain a shorter resolution that corresponds to a homomorphism
obtained from $f_i$ by bending (the torus $F_1\cup F'$ gives us the
element which determines the bending).

If $\Lambda$ has more than one family of parallel leaves, we shorten
the resolution as above working on each family at the same
time. Hence, we conclude that if $\Lambda$ is 
completely simplicial, the subgroups carried by generic leaves are all
abelian, and $f_i$'s do not factor through the quotients we made in
this case, then the
resolutions $\Phi_i$ were not shortest as assumed. 
\end{proof}

\begin{remark}
We have assumed that $\Lambda$ induces a splitting of $G$ with two
vertex groups. We realize that that assumption is not critical, as the
case of an HNN extension would be proved in the exactly same way.  It
is worth noting that bending in the case of HNN extensions, which we
perform in order to shorten the resolutions, will yield a
multiplication by elements of parabolic subgroups. This is exactly
what we needed in Example \ref{badex}.
\end{remark}

\begin{lemma}\label{surfthin}
If $\Lambda$ has a minimal component of either surface or thin type we
make quotients by trivializing all the loops in the prechosen
leaf and by abelianizing the fundamental group of the regular
neighborhood of $\Lambda$. These quotients are proper quotients of $G$. 
\end{lemma}

\begin{proof}
If $\Lambda$ has a component of thin type, then the quotient as above
has to be proper for otherwise $G$ was freely decomposable. 

Suppose now $\Lambda$
has a surface component. If the annuli that we quotiented out were
nontrivial in $G$, we obviously get a proper quotient. Suppose that
was not the case. We have a certain number of annuli attached to the
surface in our complex $K$, and to those annuli are maybe attached
different components of $K$, which we will call black boxes. Let us
remind ourselves that the annulus $S^1 \times I$ with a lamination
$S^1 \times \{\text{Cantor set}\}$ is attached along an arc transverse
to the lamination, Figure~\ref{surface}.
\begin{figure}[!ht]
\centerline{\input{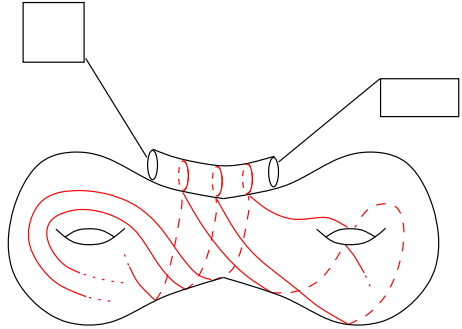}}
\caption{Lamination on the surface with the attached annuli}\label{surface}
\end{figure}

Since all the loops in the leaves of the lamination that are contained
in the annulus, i.e., $S^1 \times \{ \text{pt}\}$, are trivial we may
collapse all of them, including the circles that are not contained in
the leaves of the lamination. After this collapsing the black boxes
are attached either directly to the surface, that is to the components
complementary to the lamination, or to a boundary component of the
surface. Since the lamination on the surface is the filling one, the
components complementary to the lamination are either simply connected
or homotopy equivalent to a boundary component. Hence if the
fundamental group of one of the black boxes attached to the simply
connected complementary component is nontrivial we get a free
decomposition of $G$, which is a contradiction. We conclude that all
black boxes are either attached to the boundary components of the
surface or have trivial fundamental groups (and hence can be
collapsed).  In both cases we get a pure surface component and a surface
vertex group in a cyclic decomposition of $G$. This surface vertex group must be conjugate to a QH vertex group in the JSJ
of $G$, see \cite{mnotes}. The fact that $\Lambda$ was obtained as a limit of
$\Lambda_{f_i}$'s and that it is filling on the surface implies that
the loops that correspond to each generator of our surface group have
positive measures. Proposition 5.8. in \cite{stable} guarantees the
existence of the set of short generators. That is: by applying Process
II to the component of $K$ that carries $\Lambda$, we can obtain from
$K$ a new band complex $K'$ which resolves the same tree, but so that
the fundamental group of the component carrying the lamination is
generated by loops short with respect to an interval $s$ of
arbitrarily small measure. That the loop $p$ is short with respect to
$s$ means that $p=p_1*\lambda*p_2$, where $p_1$ and $p_2$ are
contained in $s$ and $\lambda$ is contained in a leaf of $\Lambda$.
This now tells us that these generators are shorter with respect to
the measure $\mu_i$ than the original generators were. 
Since the combinatorial type of band complexes that correspond to
these two sets of generators are the same we can find a modular
automorphism that takes one of these sets into the other.
Precomposing the $f_i$'s by this modular automorphism will give us
homomorphisms shorter than $f_i$'s, which is a contradiction with our
choice of $f_i$'s.

Finally, if $\Lambda$ has either surface or thin type component the
the fundamental group of the regular neighborhood of that component is
not abelian, hence its abelianization leads to a proper quotient of
$G$. 
\end{proof}

\begin{lemma}\label{toral}
Suppose $\Lambda$ has a minimal component $\Lambda_0$ which
resolves a line. Let $N$ be a standard neighborhood of
    $\Lambda_0$. We form the following quotients:
\vspace{-.2cm}
\begin{itemize}
    \item Abelianize $\pi_1(N)$, where $N$ is a standard neighborhood of
    $\Lambda_0$, if it is nonabelian.
\vspace{-.2cm}      
\item If $\pi_1(N)$ is abelian, make the following quotients:
        \begin{itemize}
\vspace{-.2cm}
          \item abelianize the subgroup of $G$ generated by $\pi_1(N)$
          and all vertex groups adjacent to it,
\vspace{-.cm}          
\item mod out by the direct summand of $\pi_1(N)$ that
          intersects trivially the peripheral subgroup, if such exists.
          \end{itemize}
    \end{itemize}
At least one of these quotients is proper. 
\end{lemma}

\begin{proof}

If $H=\pi_1(N)$ is nonabelian, then its abelianization leads to a
proper quotient. Assume then that $H$ is an abelian group and $N$ is a genuine
torus. If $\Lambda$ also has a simplicial, surface or thin component,
we apply Lemmas \ref{simplicial}, \ref{surfthin}. Hence our only concern is
when $\Lambda$ has only toral components, each of which has a torus as
a standard neighborhood. Since $G$ is a limit group, these tori can
not be glued directly to each other, i.e., an edge space cannot embed
into both of them.

Let $H'$ denote the peripheral subgroup of $H$. Suppose first $H'$ is
contained in a proper direct summand of $H$. We will call the
smallest such $H'$ so that $H=H' \oplus H''$.
For each $i$ we define a homomorphism $f_i':G \to L$ that coincides
with $f_i$ on all vertex groups in the decomposition of $G$ except for
$H$ where we define it to be:
\[f_i'(h)=\left\{
\begin{array}{cl}
f_i(h), & h \in H'\\
1, & h \in H''
\end{array} \right.
\]  
$f_i$ and $f_i'$ are equivalent under our moves, and $f_i'$'s factor
through the quotient $G/\ll H'' \gg$. If, on the other hand, $H'$ is a
finite index subgroup of $H$ there exists an edge $a$ of the
triangulation of $K$ intersecting an edge space $E$ adjacent to $N$
nontrivially so that $\mu_i(a) \sim d_i$. Let $V$ be the vertex space
at the other end of $E$.  Since $f_i(H)<A$ is abelian,
$\Phi_i(\widetilde{N})$ is contained in a bounded neighborhood of the
horosphere $P$ bounding the horoball $B$ whose stabilizer is $A$. We
first assume that the minimum distance between points in
$\Phi_i(\widetilde{N})$ and points in $P$ is proportional to $d_i$. 
We can find a closed path $p$ contained entirely in $V$ and $E$ whose
measure is 'virtually' 0 and whose image does not belong to the
centralizer of $f_i(H)$. If such a
path did not exist we would get a proper quotient of $G$ by modding
out the relations that make $\pi_1(V)$ and $H$ commute. We consider
the lift $\tilde{p}$ of $p$ having a nonempty intersection with
$\widetilde{N}$. One of the endpoints of $\Phi_i(\tilde{p})$ will
belong to $\Phi_i(\widetilde{N})$, but $\Phi_i(\tilde{p})$ is not
completely contained in $B$, since $f_i(p)$ does not commute with
$f_i(H)$. Hence, $\ell(\Phi_i(\tilde{p})) \geq 2d_i$, contradicting
the assumption that $\mu_i(p) \sim 0$. We know now that
\[\Phi_i(\widetilde{N}) \subset N_{n_i}(P) \ \ \ \text{and} \ \ \ \lim_{i \to
  \infty} \frac{n_i}{d_i}=0. \] 
All the hypotheses of Lemma \ref{simplicial}(1) are satisfied by our
resolutions. Therefore, we find the homomorphisms shorter than $f_i$
belonging to the same equivalence classes.

\end{proof}
\begin{lemma}
Suppose $\Lambda \in LIM'(K)$ is a limit of $\Lambda_{f_i}$, where
$f_i:G\to L$ is a sequence of short homomorphisms for which $d_{f_i}
\to \infty$. There is a neighborhood $U\subset LIM'(K)$ of $\Lambda$
such that if $\Lambda_f \in U$ then $f$ factors through one of the
quotients defined in Lemmas \ref{simplicial}, \ref{surfthin}, and \ref{toral}.
\end{lemma}

\begin{proof}

It can be shown that if a subgroup of $G$ fixes an arc in the
limiting tree $T_{\infty}$, then the image of that subgroup under
almost all $f_i$ is abelian. The proof in \cite{bs} needs only small
addition that deals with the existence of cups. 

- $\Lambda$ {\it is entirely simplicial}: 

The subgroup $H$ carried by a generic leaf fixes an arc in
$T_{\infty}$. If $H$ itself is abelian, we noted in Lemma \ref{simplicial}
that $f_i$ would factor through the quotients we made in that case. If
$H$ is not abelian, then $f_i(H)$ is, hence $f_i$ will factor through
the abelianization of $H$.

- $\Lambda$ has {\it toral components}

Any two bands will commute, and so every commutator fixes an arc in
$T_{\infty}$, i.e., $[H,H]$ fixes this arc. By the same argument as in
the simplicial case we get that, for a sufficiently large $i$,
$f_i([H,H])$ is an abelian group.  By the Tits alternative the
subgroup $f_{i}(H)$ is either virtually abelian or contains a free
group on two generators. If the latter is the case, then
$[f_{i}(H),f_{i}(H)]$ could not be an abelian group, hence $f_{i}(H)$
is virtually abelian. Now we have $f_{i}(H)$ as a virtually abelian
subgroup of a torsion-free limit group, which means it is abelian. We
conclude that $f_{i}([a,b])=1, \ \ \forall a, b \in H$.

- $\Lambda$ has a {\it surface or thin component}

The following argument is due to M. Bestvina.  We consider a standard
neighborhood $N$. A loop in an annulus will fix an arc $I$ in
$T_{\infty}$. Choose the segment $I$ so that every subsegment has the
same stabilizer, and let $g \in G$ be the element corresponding to the
loop that fixes $I$. We will argue that if $f_i(g) \neq 1$ for almost
all $i$, then almost all $f_i$ map $\pi_1(N)$ into an abelian subgroup
of $L$ and hence they all factor through an abelianization of
$\pi_1(N)$.  Suppose there is an element $h \in G$ such that $J=h(I)
\cap I$ is a nonempty interval. Then both $g$ and $g'=hgh^{-1}$ fix
$J$, and hence $I$ since the action on $T_{\infty}$ is stable. Our
remark at the beginning says that $[f_n(g),f_n(g')]=1$ for almost all
$n$. Since $f_n(h)$ conjugates $f_n(g)$ into $f_n(g')$ and maximal
abelian subgroups of $L$ are malnormal, we conclude that all three
images must belong to the same abelian subgroup of $L$. Therefore,
$[f_n(h),f_n(g)]=1$, for almost all $n$. We get that $f_n(g)$ commutes
with $f_n(h)$ whenever $h(I)\cap I\neq \emptyset$. Since $I$ belongs
to the surface component such elements $h$ generate the whole surface
(existence of ``small'' finite generating set of $\pi_1(N)$ (see
\cite{stable} -- this only uses minimality), and so $f_i(\pi_1(N))$ is
abelian for almost all $i$.

\end{proof}

\begin{proof}[Proof of Main Theorem]

Consider a sequence of homomorphisms for
which the sequence $\{d_i\}$ is bounded. Proposition \ref{finite} tells us
that in such a sequence we can have only finitely many nonconjugate
homomorphisms. We pick a representative of each conjugacy class, say
$k_i:G\to L, \ i=1, \ldots, k$ and form quotients $K_i=G/ker(k_i)$. 
If this quotient is not proper, then $k_i$ was an embedding.

If on the other hand, $\Lambda \in LIM'(K)$ such that
$\Lambda_{f_i}\to \Lambda$ and $d_{f_i}\to \infty$, then our previous
four lemmas prove the claim: we have formed finitely many quotients
of $G$ and we have found a neighborhood $U_{\Lambda}$ of $\Lambda$ so
that whenever $\Lambda_f \in U_{\Lambda}$ and $d_{[f]}\gg 0$ then a
homomorphism equivalent to $f$ factors through one of these
quotients. $U_{\Lambda}$'s together with neighborhoods of
$\Lambda_{k_i}$'s cover $LIM'(K)$. Since this space is compact, it is
covered by finitely many of these neighborhoods. Hence we have finitely
many quotients through which an element of the equivalence class
$\sim$ of any $f:G \to L$ with $d_f>0$ factors. If $d_f=0$, then $f$
factors through abelianization of $G$.

This concludes the proof of our theorem. 
\end{proof}

It would be useful to know what happens if we iterate
this construction, i.e., if we apply Theorem \ref{mainthm} to all the
quotients we obtained. The following lemma is well known. 

\begin{lemma}\label{epi}
A sequence of epimorphisms between $\omega$-residually free groups
eventually stabilizes.
\end{lemma}

\begin{proof}

Let
\[G_1 \to G_2 \to \cdots \to G_n \to \cdots \]
be a sequence of epimorphisms between $\omega$-residually free
groups. For a fixed free group $F$ the sequence 
\[Hom(G_1,F) \leftarrow Hom(G_2,F)\leftarrow \cdots
Hom(G_n,F)\leftarrow \cdots \] eventually stabilizes, i.e., consists of
bijections, see \cite{mnotes}. Suppose 
\[Hom(G_k,F) \cong Hom(G_{k+1},F), \]
where the isomorphism is given by the obvious inclusion of
$Hom(G_{k+1},F)$ into $ Hom(G_k, F)$. Let us also suppose that
$e: G_k \to G_{k+1}$ is a proper epimorphism. Let $g \in
ker(e)$. Since $G_k$ is $\omega$-residually free there is a
homomorphisms $f:G_k \to F$ such that $f(g) \neq1$. On the other hand,
there exists $f':G_{k+1} \to F$ such that $f=f' \circ e$. We now have
\[1\neq f(g)=f'(e(g))=1. \]
Hence, $G_k \cong G_{k+1}$. 
 
\end{proof}

We now form the Makanin-Razborov diagram, Figure \ref{mr}, except
we add an edge issuing from each group in this diagram and ending in
$L$ representing the embeddings. All the groups in this diagram,
except possibly $G$, will be limit groups. Furthermore, each branch of
this diagram is finite, as we have just shown.

\begin{remark}
At a first glance it may appear as if the proof of Main Theorem should
hold for a broader class of groups, in particular for groups that are 
hyperbolic relative to a collection of their maximal noncyclic abelian
subgroups. There are two problems that appear. The first one is that
we do not know that a finitely generated subgroup of a relatively
hyperbolic group is finitely presented, and so the proof would need to
be modified in order to deal with not necessarily finite
complexes. The second problem is that the proof of Lemma \ref{epi}
does not apply in this context. However, we believe that the
proof of descending chain condition on (hyperbolic)-limit groups given
in \cite{zlil8} will transfer to this setting easily. Nonetheless,
until a more elegant solution is found to both of these problems, we
leave the discussion at the level of this remark.

\end{remark}

\bibliographystyle{siam}
\bibliography{all}
\par

\end{document}